\documentclass[12pt,a4paper,reqno]{amsart}

\newcommand{\nc}[2]{\newcommand{#1}{#2}}
\newcommand{\rnc}[2]{\renewcommand{#1}{#2}}

\nc{\ot}{\otimes}
\nc{\tot}{\hat\ot}
\nc{\opl}{\oplus}
\nc{\liea}{\mathfrak{g}}	
\nc{\sliea}{\mathfrak{h}}
\nc{\uea}{U(\liea)}
\nc{\quea}{U_{h}(\liea)}
\nc{\sln}{\mathfrak{sl}_{n}(\CC)}
\nc{\sls}{\mathfrak{sl}_{3}(\CC)}
\nc{\sun}{\mathfrak{su}(l+1)}
\nc{\qusln}{U_{h}(\sln)}
\nc{\sonn}{\mathfrak{so}_{2l+1}(\CC)}
\nc{\sonnr}{\mathfrak{so}(2l+1)}
\nc{\spn}{\mathfrak{sp}_{2l}(\CC)}
\nc{\spnr}{\mathfrak{sp}(2l)}
\nc{\son}{\mathfrak{so}_{2l}(\CC)}
\nc{\sonr}{\mathfrak{so}(2l)}
\nc{\glv}{\mathfrak{gl}(V)}
\nc{\gll}{\mathfrak{gl}(\mathfrak{g})}
\nc{\gln}{\mathfrak{gl}_{n}(\CC)}
\nc{\glt}{\mathfrak{gl}_{3}(\CC)}
\nc{\slfg}{\mathfrak{sl}_{2N-1}(\CC)}
\nc{\glf}{\mathfrak{gl}_{4}(\CC)}
\nc{\glnr}{\mathfrak{gl}(l+1)}
\nc{\CC}{\mathbb{C}}
\nc{\QQ}{\mathbb{Q}}
\nc{\RR}{\mathbb{R}}
\nc{\QQbar}{\overline{\QQ}}
\nc{\ur}{\mathcal{R}}
\nc{\uf}{\mathcal{F}}
\rnc{\d}{\Delta}
\rnc{\b}{\beta}
\nc{\s}{\sigma}
\rnc{\l}{\lambda}
\nc{\ep}{\epsilon}
\nc{\fa}{\forall}
\nc{\lah}{\varphi}
\nc{\lgh}{\phi}
\nc{\mapto}{\rightarrow}
\nc{\act}{\triangleright}
\nc{\ract}{\triangleleft}
\nc{\ch}{\mathfrak{h}}
\rnc{\a}{\alpha}
\rnc{\k}{\kappa}
\nc{\llra}{\Longleftrightarrow}
\theoremstyle{plain}
\newtheorem{thm}{Theorem}
\newtheorem{prop}[thm]{Proposition}
\newtheorem*{lem}{Lemma}

\theoremstyle{definition}
\newtheorem{defn}[thm]{Definition}
\newtheorem{ex}[thm]{Example}

\theoremstyle{remark}
\newtheorem{remark}[thm]{Remark}

\numberwithin{thm}{section}

\nc{\bth}{\begin{thm}}
\nc{\bprop}{\begin{prop}}
\nc{\ble}{\begin{lem}}
\nc{\bde}{\begin{defn}}
\nc{\bex}{\begin{ex}}
\nc{\bre}{\begin{remark}}

\nc{\ethe}{\end{thm}}
\nc{\eprop}{\end{prop}}
\nc{\ele}{\end{lem}}
\nc{\ede}{\end{defn}}
\nc{\eex}{\end{ex}}
\nc{\ere}{\end{remark}}

\DeclareMathOperator{\op}{op}

\DeclareMathOperator{\id}{id}

\DeclareMathOperator{\Hom}{Hom}
\DeclareMathOperator{\Img}{Im}

\title[Twisting $2$-cocycles]{Twisting $2$-cocycles for the construction of 
new non-standard quantum groups}
\author{Andrew D. Jacobs and J. F. Cornwell}
\address{Department of Physics and Astronomy, University of St.~Andrews,
North Haugh, St.~Andrews, Fife, KY16 9SS, Scotland}
\begin{document}

\begin{abstract}
We introduce a new class of $2$-cocycles defined explicitly on the generators of
certain 
multiparameter standard quantum groups. These allow us, through the process of 
twisting the familiar standard quantum groups,
to generate new as well as previously known examples of non-standard 
quantum groups. In particular we are able to construct generalisations 
of both the Cremmer-Gervais deformation of $SL(3)$ and 
the so called esoteric quantum groups
of Fronsdal and Galindo in an explicit and straightforward manner.
\end{abstract}

\maketitle

\section{Introduction}

Originally there were two clearly defined types of quantum 
groups~\cite{D0,J,FRT}. They were \emph{single-parameter} quantisations, 
$U_{q}(\liea)$ and $\CC_{q}[G]$ respectively, of dual classical objects: the 
universal enveloping algebras of simple Lie algebras, $U(\liea)$, and the 
coordinate rings of simple Lie groups, $\CC[G]$. With their universal $R$-matrices, 
$\mathcal{R}$, these $U_{q}(\liea)$ are the standard examples of quasitriangular 
Hopf algebras, while the $\CC_{q}[G]$, together with the corresponding 
numerical $R$-matrices, are the standard examples of what we call 
co-quasitriangular Hopf algebras~\cite{D,DT}. It soon became 
apparent that there were a number of multiparameter 
generalisations~\cite{DMMZ,Sud,Tak} of these standard quantum groups and 
through the work of Drinfeld~\cite{D1}, followed by Reshetikhin~\cite{Resh}, 
an interpretation emerged: all the multiparameter quantum groups corresponding 
to a particular standard quantum group were related, amongst themselves 
and with the standard quantum group, through Drinfeld's important process of 
twisting. In fact the original works of Drinfeld and Reshetikhin were 
concerned only with quasitriangular Hopf algebras, but their constructions 
dualise immediately to the case of co-quasitriangular Hopf algebras. 
Since the twists act only as similarity transformations on the so called 
$\hat{R}$-matrices~\cite{FRT}, the different standard quantum groups 
corresponding to different classical Lie groups cannot be related to each 
other by twisting. The picture then is of a number of distinct 
`twist equivalence classes'. Later, Kempf and 
Engeldinger~\cite{K,Eng} (see also the work of Khoroshkin and Tolstoy~\cite{KT}) refined 
Reshetikhin's work slightly and showed that there were other 
interesting quantum groups related, through Reshetikhin-type twists, 
with the standard ones.

From time to time there appeared genuinely non-standard quantum 
groups, usually defined in terms of non-standard numerical $R$-matrices. 
It is natural to investigate whether these define new twist equivalence classes 
or whether they belong to classes already defined by the standard $R$-matrices.  
We will be particularly concerned in this article with the non-standard quantum 
groups of Cremmer and Gervais~\cite{CG} and  
Fronsdal and Galindo~\cite{FG1,FG2} which for general theoretical reasons 
(see Section~4) 
may be expected to be twist-equivalent to the standard $SL(n)$ 
quantum groups. However, let us be clear that \emph{no} twist of the 
Reshetikhin type is suitable in these cases. As we explain later, 
the relevant twisting 
structures are counital 2-cocycles on the quantum groups.
The problem then is to find the appropriate twisting 2-cocycles, 
defined on the standard co-quasitriangular Hopf algebra, $\CC_{q}[SL(n)]$, 
which twist this quantum group into the Cremmer-Gervais and Fronsdal-Galindo quantum groups. 

Recently, Hodges has made significant progress in this area~\cite{H4}, drawing on 
previous work of his contained in a series of important 
papers~\cite{H1,H2,H3}. These covered many aspects of quantum group theory 
from a ring theoretic perspective. We should mention~\cite{H2} in particular,
where some remarkable aspects of the algebraic structure of Cremmer-Gervais 
quantum groups were revealed. 
In~\cite{H4} Hodges starts from a particular, standard, multiparameter
quantized enveloping algebra $U_p(\liea)$. He then identifies a pair of
\emph{commuting} sub-Hopf algebras, $U_{p}(\mathfrak{b}_{1}^{-})$ and
$U_{p}(\mathfrak{b}_{2}^{+})$, associated with certain Belavin-Drinfeld 
triples~\cite{BD}. This gives rise to a Hopf algebra homomorphism,
$\phi:U_{p}(\mathfrak{b}_{1}^{-})\ot U_{p}(\mathfrak{b}_{2}^{+})\mapto U_p(\liea)$,
through the usual multiplication map. Attention then shifts to the dual
map, $\phi^{*}:\CC_p[G]\mapto \CC_{p}[B_{1}^{-}]\ot \CC_{p}[B_{2}^{+}]$.
Hodges proceeds to identify $\Img(\phi^{*})$, in a series of precise and 
subtle steps, with the image of the tensor product of 
a pair of `extended' Borel subalgebra-like objects, between which
there is a skew pairing. This skew pairing lifts to the tensor components
of $\Img(\phi^{*})$. It is well known that such a pairing gives rise
to a $2$-cocycle, the quintessential example appearing in the twisting
interpretation of the quantum double~\cite{DT,Majb}, and a $2$-cocycle
is then induced on the quantized function algebra $\CC_p[G]$. 
Hodges claims that in the particular case of $\sls$ the $2$-cocycle coming from 
his construction generates the Cremmer-Gervais deformation of $\CC[SL(3)]$. More
generally, he claims that it should also be possible to reach the esoteric 
quantum groups of Fronsdal and Galindo~\cite{FG1,FG2}. However, 
Hodges' approach is rather technical and does not readily yield 
$2$-cocycles defined \emph{explicitly} on the familiar $T$ generators of 
standard quantum groups. We are able to remedy this situation here.

Our approach is actually quite distinct from that of Hodges. 
We work entirely within the framework of co-quasitriangular Hopf algebras 
coming from solutions 
of the matrix quantum Yang-Baxter equation (QYBE). In Section~2 we recall the 
definition of a co-quasitriangular Hopf algebra and
the basic result on twisting by $2$-cocycles. A good reference for this theory,
and much more besides, is the book by Majid~\cite{Majb}, from which much of 
our notation is borrowed. Other good references
are the paper by Larson and Towber~\cite{LT} and the papers
of Doi and Takeuchi~\cite{D,DT}. Majid calls co-quasitriangular Hopf algebras, 
`dual quasitriangular Hopf algebras', but our terminology comes 
from~\cite{D,DT}. We go on to describe the well
known class of 2-cocycles which appears as the dual of the 
Reshetikhin-type twists. They originate from particular solutions of the
QYBE.  The `parameterization' twists originally considered by Reshetikhin~\cite{Resh}
may be regarded as examples in this class, and using such a twist we have 
obtained a \emph{new} $3$-parameter generalized Cremmer-Gervais $R$-matrix, presented 
here, which includes as a special case the $2$-parameter $R$-matrix considered by 
Hodges in~\cite{H2}. Details of the derivation of this new $R$-matrix are 
given in Appendix~A. The sub-Hopf-algebra-induced twists considered
by Engeldinger and Kempf~\cite{K,Eng} also belong to this general class
of $2$-cocycles, and are recalled here. 

Section~3 contains our main new results. We present there a new class of
$2$-cocycles which no longer emanate from solutions of the matrix QYBE. 
Instead, they arise from matrices satisfying a new,
and remarkably simple, system of equations. A number of explicit 
$2$-cocycles belonging to our new class are presented, along with
the following results:
\begin{itemize}
\item
It is shown \emph{explicitly} that
the \emph{new} $3$-parameter generalised Cremmer-Gervais quantum group corresponding 
to $GL(3)$, already given in Section~2, is obtained from a particular 
multiparameter standard quantum group through twisting. 
\item
The $2$-cocycle used to obtain the generalised Cremmer-Gervais deformation of $GL(3)$
is an example of a general class of \emph{simple root} $2$-cocycles
which themselves belong to a more general class of 
\emph{composite simple root} $2$-cocycles. These $2$-cocycles may
all be defined on certain standard, multiparameter, deformations of
$GL(n)$, and consequently generate new non-standard quantum groups.
\item
A $2$-cocycle which can be used to twist a certain multiparameter
standard deformation of $\CC[GL(2N-1)]$ to obtain a \emph{generalisation}
of the quantum groups of Fronsdal and Galindo is presented. For $N=2$ 
this $2$-cocycle is just the one used to obtain the generalised Cremmer-Gervais 
$GL(3)$ quantum group.
\end{itemize}
Note that the $R$-matrix considered by Fronsdal and Galindo is already
`multiparameter', involving $N$ parameters, but we obtain, in Appendix~B,
an $R$-matrix which , for $N>3$, depends on $(1+\frac{1}{2}(N-1)(N+4))$ 
parameters and, in the cases $N=2$ and $N=3$, on $3$ and $7$ 
parameters respectively. Let us also note that, as was already suggested in Hodges 
work~\cite{H4}, starting from the original standard quantum groups, $\CC_{q}[SL(n)]$, 
we need a \emph{combination} of the Reshetikhin-type parameterisation twists with 
our new twists to obtain the Fronsdal-Galindo quantum groups.

In Section~4, we collect some information about the semi-classical objects
corresponding to the $R$-matrices which we have been considering in this paper,
namely the classical $r$-matrices.
We also recall the background, in Drinfeld's fundamental work, which
serves as the on-going motivation in the quest for interesting twists. We end
by pointing out a particular problem involved in constructing the 
Cremmer-Gervais quantum group corresponding to $GL(4)$, and describe 
an interesting phenomenon involving sub-Hopf-algebra-induced twists. Starting
from the standard multiparameter quantum group, $\CC_{q,p}[GL(4)]$, 
we twist, first of all, using a sub-Hopf-algebra-induced twist. The resulting 
quantum group may reasonably be called `weakly non-standard' and can be 
twisted further using a $2$-cocycle from our construction. The new, 
non-standard, $R$-matrix obtained through this double twist involves a pair
of non-standard off-diagonal elements which could not have been added to the 
original $R$-matrix directly, but which \emph{do} appear in the Cremmer-Gervais
$R$-matrix for $GL(4)$.

\pagebreak

\section{Co-quasitriangular Hopf algebras and Reshetikhin twists}
We begin with some basic definitions.
\bde
A bialgebra $A$ is called 
\emph{co-quasitriangular} if there exists a bilinear form $\s$ on $A$, 
which we will
call an $R$-form, such that
\begin{enumerate}
\item $\s$ is invertible with respect to the convolution product, $*$, 
that is, there is 
another bilinear form $\s^{-1}$ such that
\begin{equation}
\s(a_{(1)},b_{(1)})\s^{-1}(a_{(2)},b_{(2)})=\ep(ab)=
\s^{-1}(a_{(1)},b_{(1)})\s(a_{(2)},b_{(2)}),  \label{cq1}
\end{equation}
\item $\s*m=m^{\op}*\s$, i.e.
\begin{equation}
\s(a_{(1)},b_{(1)})a_{(2)}b_{(2)}=b_{(1)}a_{(1)}\s(a_{(2)},b_{(2)}), 
\label{cq2}
\end{equation}
\item $\s(m\ot\id)=\s_{13}*\s_{23}$, i.e.
\begin{equation}
\s(ab,c)=\s(a,c_{(1)})\s(b,c_{(2)}), 
\end{equation}
\item $\s(\id\ot\, m)=\s_{13}*\s_{12}$, i.e.
\begin{equation}
\s(a,bc)=\s(a_{(1)},c)\s(a_{(2)},b).
\end{equation}
\end{enumerate}
\ede

\bre
We are employing here a slightly simplified version of the Sweedler 
notation for coproducts: $\d(a)=a_{(1)}\ot a_{(2)}$, with the 
summation suppressed.
\ere

\bre
It may be useful to briefly recall how we arrive at this definition.
Suppose that $(H,\mathcal{R})$ is a quasitriangular bialgebra, with
	$\mathcal{R}\in H\ot H$ the \emph{universal $R$-matrix} obeying
Drinfeld's familiar axioms,
\begin{gather}
(\d\ot\id)(\mathcal{R})=\mathcal{R}_{13}\mathcal{R}_{23},\;\;\;\;\;
(\id\ot\d)(\mathcal{R})=\mathcal{R}_{13}\mathcal{R}_{12}, \\
\d^{\op}(h)=\mathcal{R}\circ\d(h)\circ\mathcal{R}^{-1},\;\;\;\fa h\in H.
\end{gather}
In fact, we may be regard $\mathcal{R}$ as a map $k\mapto H\ot H$,
where $k$ is the ground field. When formulating the dual notion of co-quasitriangular bialgebra, 
we then need to consider an $R$-form $\s$ which is now a map $A\ot A\mapto k$,
where $A$ can be thought of as dual to $H$. To formulate the dual axioms,
involving $\s$ instead of $\mathcal{R}$,
we require the algebra structure on $\Hom(A\ot A,k)$. This is the convolution
algebra provided by the natural tensor product coalgebra structure of $A\ot A$.
Explicitly then, let us give the details for a particular example,
\begin{align}
(\s_{13}*\s_{23})(a\ot b\ot c)&=\s_{13}(a_{(1)}\ot b_{(1)}\ot c_{(1)})
\s_{23}(a_{(2)}\ot b_{(2)}\ot c_{(2)}) \nonumber \\
&=\s(a_{(1)},c_{(1)})\ep(b_{(1)})\s(b_{(2)},c_{(2)})\ep(a_{(2)}) \nonumber \\
&=\s(a,c_{(1)})\s(b,c_{(2)}). \nonumber
\end{align}
For more on this process of `dualising' we refer the reader to Majid's
book~\cite{Majb}, and his paper~\cite{Maj1}.
\ere

From the definition it is readily seen that the QYBE now manifests itself as
\begin{equation}
\s_{12}*\s_{13}*\s_{23}=\s_{23}*\s_{13}*\s_{12}. \label{yb}
\end{equation}
When $A$ is actually a Hopf algebra, with
antipode $S$, we call it a \emph{co-quasitriangular Hopf
algebra}. It can then be shown that $S$ is \emph{always} 
invertible, and
$\s^{-1}(a,b)=\s(S(a),b)$ and $\s(a,b)=\s^{-1}(a,S(b))$. 

The FRT bialgebras $A(R)$, introduced by the Leningrad school~\cite{FRT} 
and developed by Majid~\cite{MajFRT}, where $R$ is any matrix solution
of the QYBE, fit into
the co-quasitriangular bialgebra framework. Indeed, we define the 
$R$-form on the generators $T_{i}^{j}$, as 
\begin{equation}
\s(T_{i}^{s},T_{j}^{t})=R_{ij}^{st},
\end{equation}
or in the useful `matrix notation', as 
\begin{equation}
\s(T_{1},T_{2})=R_{12},
\end{equation}
and then extend its domain of definition to the
whole of $A(R)$ by setting
\begin{align}
\s(T_{1}T_{2},T_{3})& =\s(T_{1},T_{3})\s(T_{2},T_{3})=R_{13}R_{23},\label{bichar1} \\
\s(T_{1},T_{2}T_{3})& =\s(T_{1},T_{3})\s(T_{1},T_{2})=R_{13}R_{12}.\label{bichar2}
\end{align}
The QYBE then guarantees consistency with the product relation (\ref{cq2}).
\bre
To remove any possible doubt about the notation being employed here, let us
present~(\ref{bichar1}) explicitly, in terms of the generators, as
\begin{equation}
\s(T_{i}^{s}T_{j}^{t},T_{k}^{r})=\s(T_{i}^{s},T_{k}^{m})\s(T_{j}^{t},T_{m}^{r})
=R_{ik}^{sm}R_{jm}^{tr}, \nonumber
\end{equation}
where the summation convention \emph{is} being assumed.
\ere

\bde
A bilinear form $\chi$ on a bialgebra $A$ is called a \emph{counital 2-cocycle} 
on $A$ if it is invertible in the convolution product, and
\begin{equation}
\chi(1,a)=\ep(a)=\chi(a,1) \label{coc1},
\end{equation}
and
\begin{equation}
\chi_{12}*\chi(m\ot\id)=\chi_{23}*\chi(\id\ot\, m). \label{coc2}
\end{equation}
\ede

\bre
It is a simple matter to show that any $R$-form is a counital 2-cocycle.
We also note that for any Hopf algebra on which we can define 
such a 2-cocycle, which moreover intertwines the multiplication
as in~(\ref{cq2}), the antipode is necessarily invertible. 
\ere

\bre
In the more familiar dual version of this definition, we consider
an invertible element $\mathcal{F}$ of $H\ot H$. Then~(\ref{coc1}) 
and~(\ref{coc2}) correspond respectively to
\begin{equation}
(\ep\ot\id)(\mathcal{F})=1=(\id\ot\ep)(\mathcal{F}), \label{dcoc1}
\end{equation}
and 
\begin{equation}
\mathcal{F}_{12}(\d\ot\id)(\mathcal{F})=\mathcal{F}_{23}(\id\ot\d)(\mathcal{F}).
\label{dcoc2}
\end{equation}
An element $\mathcal{F}$ satisfying these conditions is then called a
counital 2-cocycle \emph{for} $H$.
\ere

\bre
For a general discussion of cocycles for and on Hopf algebras we refer the 
reader to Section~2.3 of Majid's book~\cite{Majb}.
\ere

The property of co-quasitriangular Hopf algebras which is of 
particular interest to us is that, given one, we may generate others
using these counital 2-cocycles.
This important process of \emph{twisting} is the dual of Drinfeld's 
original quasitriangular quasi-Hopf algebra twist~\cite{D1},
restricted to the special case of twisting \emph{from and to} 
co-quasitriangular 
Hopf algebras. It is not difficult to dualise Drinfeld's 
original quasitriangular quasi-Hopf algebra axioms, 
and his result on twisting.
This was probably first carried out explicitly by Majid~\cite{Maj1}.
We obtain the axioms for a co-quasitriangular co-quasi-Hopf algebra and on 
specialising the twisting result, we obtain the following important theorem.
\bth
Let $(A,m,\eta,\d,\ep,\s)$ be a co-quasitriangular bialgebra and let
$\chi$ be a counital 2-cocycle on $A$, then there is a new 
co-quasitriangular bialgebra $(A_{\chi},\s_{\chi})$ obtained by
twisting the product and $R$-form of $(A,\s)$ as
\begin{align}
m_{\chi}& =\chi*m*\chi^{-1}, \label{twstm}\\
\s_{\chi}& =\chi_{21}*\s*\chi^{-1}. \label{twstr}
\end{align}
If $A$ is moreover a Hopf algebra with antipode $S$, then $A_{\chi}$ is
also a Hopf algebra with twisted antipode given by
\begin{equation}
S_{\chi}=\l*S*\l^{-1}, \label{twsta}
\end{equation}
where $\l=\chi\circ(\id\ot\, S)\circ\d$.
\ethe

For the co-quasitriangular bialgebras, $A(R)$, there is a particularly obvious 
way of constructing twisting 2-cocycles. 
Take any invertible solution $F$ of the QYBE and define a bilinear form $\chi$
by 
\begin{gather}
\chi(T_{1},T_{2})=F_{12}, \\
\chi(1,T)=\chi(T,1)=\ep(T), \\
\intertext{and}
\chi^{-1}(T_{1},T_{2})=F_{12}^{-1}. 
\end{gather}
We then extend this to the whole of $A(R)$ just as we did for the $R$-form
in equations~(\ref{bichar1}) and~(\ref{bichar2}), that is
\begin{align}
\chi(T_{1}T_{2},T_{3})& =\chi(T_{1},T_{3})\chi(T_{2},T_{3})=F_{13}F_{23},  \label{cbc1} \\
\chi(T_{1},T_{2}T_{3})& =\chi(T_{1},T_{3})\chi(T_{1},T_{2})=F_{13}F_{12}. \label{cbc2}
\end{align}
However $\chi$ must respect the 
algebra structure already on $A(R)$ so we must also have
\begin{align}
\;& \chi(R_{12}T_{1}T_{2}-T_{2}T_{1}R_{12},T_{3})=0 \nonumber \\ 
\nonumber \\
\llra\;\; & R_{12}F_{13}F_{23}=F_{23}F_{13}R_{12}, \label{compat1} \\
\intertext{and}
\;& \chi(T_{1},R_{23}T_{2}T_{3}-T_{3}T_{2}R_{23})=0 \nonumber  \\ 
\nonumber \\
\llra\;\; & R_{23}F_{13}F_{12}=F_{12}F_{13}R_{23}. \label{compat2}
\end{align}
Thus any invertible solution $F$ of the QYBE which 
satisfies~(\ref{compat1})~and~(\ref{compat2}) provides a 2-cocycle twist. 
Such a twisting system $(R,F)$ may reasonably be called a 
\emph{Reshetikhin twist}~\cite{Resh}.

\bre
In the context of the papers~\cite{Resh,K,Eng}, where the approach
is dual to ours, an invertible element $\mathcal{F}\in H\ot H$ is
considered, where $(H,\mathcal{R})$ is a quasitriangular Hopf algebra.
It is assumed to satisfy the QYBE, 
\begin{equation}
\mathcal{F}_{12}\mathcal{F}_{13}\mathcal{F}_{23}=
\mathcal{F}_{23}\mathcal{F}_{13}\mathcal{F}_{12},
\end{equation}
and the relations 
\begin{align}
(\d\ot\id)(\mathcal{F})&=\mathcal{F}_{13}\mathcal{F}_{23},\\
\intertext{and}
(\id\ot\d)(\mathcal{F})&=\mathcal{F}_{13}\mathcal{F}_{12},
\end{align}
which correspond respectively to equations~(\ref{cbc1}) and~(\ref{cbc2}). 
This $\mathcal{F}$ is then a $2$-cocycle for $H$ in the sense
of Remark~2.7, and twists the \emph{comultiplication}, 
\emph{universal $R$-matrix} and antipode as
\begin{align}
\d_{\mathcal{F}}(h)&=\mathcal{F}\d(h)\mathcal{F}^{-1},
\;\;\;\;\fa h\in H, \label{tc}\\
\mathcal{R}_{\mathcal{F}}&=\mathcal{F}_{21}\mathcal{R}\mathcal{F}^{-1},
\label{tr}\\
\intertext{and}
S_{\mathcal{F}}(h)&=vS(h)v^{-1},\;\;\;\;\fa h\in H,\label{ta}
\end{align}
where $v=m\circ(\id\ot S)(\mathcal{F})$. This is actually a slight generalisation of
the presentation of Reshetikhin~\cite{Resh}, due to Kempf and Engeldinger~\cite{K,Eng}.
\ere

We will be particularly interested in the situation pertaining when we
take $R$ to be the standard $SL(n)$ type $R$-matrix given by,
\begin{equation}
(R_{S})_{ij}^{st}=\begin{cases}
q& \text{$i=j=s=t$},\\
1& \text{$i=s\neq j=t$},\\
(q-q^{-1})& \text{$i=t<j=s$}.  \label{rs}
\end{cases}
\end{equation}
In this case, with the identification of the central quantum determinant,
$A(R)$ becomes a Hopf algebra. Indeed it is the standard 
quantization of the coordinate ring of the Lie group $SL(n)$, 
denoted $\CC_q[SL(n)]$. 

Let us also present here the expressions, in our notation,
for the Cremmer-Gervais $SL(n)$ $R$-matrix and the Fronsdal-Galindo
$GL(2N-1)$ $R$-matrix. First, the Cremmer-Gervais $R$-matrix will be taken to
be 
\begin{equation}
(R_{CG})_{ij}^{st}=\begin{cases}
q& \text{$i=j=s=t$},\\
qq^{-2(j-s)/n}& \text{$i=s<j=t$},\\
q^{-1}q^{-2(j-s)/n}& \text{$i=s>j=t$}, \\
(q-q^{-1})& \text{$i=t<j=s$},\\
(q-q^{-1})q^{-2(j-s)/n}& \text{$i<s<j$, and $t=i+j-s$},\\
-(q-q^{-1})q^{-2(j-s)/n}& \text{$j<s<i$, and $t=i+j-s$}. \label{cg}
\end{cases}
\end{equation}
If the $R$-matrix which appears as equation (46) in the original paper
of Cremmer and Gervais~\cite{CG} is denoted $\tilde{R}$ then 
$R_{CG}=\tilde{R}_{21}$ (with $e^{-ih}$ there replaced be $q$ here).
$A(R_{CG})$ again becomes a Hopf algebra, with the identification
of the central quantum determinant found in~\cite{H2}, and
will then be denoted $\CC_{CG,q}[SL(n)]$.

The Fronsdal-Galindo $R$-matrix will be considered in the next section.
However we will present it here so that the reader might easily compare 
it with the Cremmer-Gervais $R$-matrix. Thus, we take the Fronsdal-Galindo 
$R$-matrix to be
\begin{equation}
(R_{FG})_{ij}^{st}=\begin{cases}
q& \text{$i=j=s=t$},\\
q& \text{$i=s=2N-j$, $j=t$, $0<j<N$},\\
q^{-1}& \text{$i=s$, $j=t=2N-i$, $0<i<N$},\\
1& \text{$i=s$, $j=t$, $i\neq j$, $i+j\neq 2N$},\\
q-q^{-1}& \text{$i=t<j=s$},\\
q\k_{i}& \text{$0<i<N$, $j=2N-i$, $s=t=N$},\\
q\tilde{\k}_{j}& \text{$0<j<N$, $i=2N-j$, $s=t=N$},\\
q^{-1}\xi_{is}& \text{$0<i<s<N$, $j=2N-i$, $t=2N-s$},\\
q\tilde{\xi}_{jt}& \text{$0<j<t<N$, $i=2N-j$, $s=2N-t$},\\
\end{cases} \label{FGR}
\end{equation}
where 
\begin{align}
\tilde{\k}_{i}&=-q^{2(N-i)}\k_{i},\\
\tilde{\xi}_{ij}&=(1-q^{-2})q^{2(j-i)}(\k_{i}/\k_{j}),\\
\xi_{ij}&=(1-q^{2})(\k_{i}/\k_{j}).
\end{align}
Clearly, $R_{FG}$ depends on $N$ parameters --- $q$ together with 
$\k_{i}$, $0<i<N$. In this case, if the $R$-matrix given 
in section 5 of the paper~\cite{FG1} (with the $q$ there replaced by $q^{-1}$)
is denoted $\tilde{R}$, then $R_{FG}=q\tilde{R}_{21}$.  It will be shown
that this $R$-matrix is related via twisting to $R_{S}$. Thus we can
say that $A(R_{FG})$ may also be taken to be a Hopf algebra, which is
more than is claimed in~\cite{FG1,FG2}. We will denote this quantum group
by $\CC_{FG,q}[GL(2N-1)]$.

\bex
For the standard deformation $\CC_{q}[SL(n)]$, defined by the $R$-form
$\s_{S}(T_{i}^{s},T_{j}^{t})=(R_{S})_{ij}^{st}$, we define the $2$-cocycle $\chi$
by
\begin{equation}
\chi(T_{i}^{s},T_{j}^{t})=F_{ij}^{st}=f_{ij}\delta_{i}^{s}\delta_{j}^{t}, \label{mpc}
\end{equation}
extended to the whole of 
$\CC_{q}[SL(n)]$ by~(\ref{cbc1}) and~(\ref{cbc2}). Note that the
summation convention is not being assumed here and indeed will not be assumed 
anywhere, unless stated otherwise.
We find that the conditions~(\ref{compat1}) and~(\ref{compat2}) impose
no restrictions on the $f_{ij}$s. The new, twisted $R$-form then coming 
from~(\ref{twstr}) is
\begin{equation}
\s_{S,p}(T_{i}^{s},T_{j}^{t})=(R_{S,p})_{ij}^{st},
\end{equation}
where $R_{S,p}$ is the familiar $(1+\binom{n}{2})$-parameter standard $R$-matrix 
given by
\begin{equation}
(R_{S,p})_{ij}^{st}=\begin{cases}
q& \text{$i=j=s=t$},\\
p_{ij}& \text{$i=s\neq j=t$},\\
(q-q^{-1})& \text{$i=t<j=s$},  \label{mprs}
\end{cases}
\end{equation}
with $p_{ij}=f_{ji}f_{ij}^{-1}=p_{ji}^{-1}$ for $i<j$. The multiparameter Hopf
algebra so defined will be denoted $\CC_{q,p}[GL(n)]$, and was first
constructed in this way by Kempf in~\cite{K} (see also the paper by
Schirrmacher~\cite{Sch}). We will often take
$p_{ii}=q$ in what follows.
\eex

\bex
\label{mpcgex}
For the $1$-parameter Cremmer-Gervais deformation $\CC_{CG,q}[SL(n)]$ 
the situation is more interesting. We again define a $2$-cocycle $\chi$
in terms of a diagonal matrix $F_{ij}^{st}=f_{ij}\delta_{i}^{s}\delta_{j}^{t}$.
However now the compatibility conditions~(\ref{compat1})~and~(\ref{compat2}) 
\emph{do} impose restrictions on the $f_{ij}$s. The number of independent parameters appearing
in the twisted Hopf algebra $\CC_{CG,q,p}[GL(n)]$ is then determined by the number
of independent combinations of the $f_{ij}$s which appear in 
$R_{CG,p}=F_{21}R_{CG}F^{-1}$. As demonstrated in Appendix~A,
we are left with just
three independent parameters --- $q$ together with a new pair, $p$ and $\l$. 
Explicitly, the \emph{$3$-parameter
generalised Cremmer-Gervais $R$-matrix} is given by
\begin{equation}
(R_{CG,p})_{ij}^{st}=\begin{cases}
q& \text{$i=j=s=t$},\\
p^{j-s}q& \text{$i=s<j=t$},\\
p^{j-s}q^{-1}& \text{$i=s>j=t$}, \\
(q-q^{-1})& \text{$i=t<j=s$},\\
p^{j-s}\l^{st-ij}(q-q^{-1})& \text{$i<s<j$, and $t=i+j-s$},\\
-p^{j-s}\l^{st-ij}(q-q^{-1})& \text{$j<s<i$, and $t=i+j-s$}. \label{mpcg}
\end{cases}
\end{equation}
We refer the reader to Appendix~A for the proof of this result.
\eex

\bex
Another type of Reshetikhin twist, is the sub-Hopf-algebra-induced twist, 
studied in particular by Engeldinger and Kempf~\cite{Eng}. An example of such
a twist is given by defining a $2$-cocycle
on $\CC_{q,p}[GL(n)]$ as 
\begin{equation}
\chi(T_{i}^{s},T_{j}^{t})=\begin{cases}
f_{ij}& \text{$i=s$, $j=t$},\\
q^{-1}(q-q^{-1})f_{\eta\eta}& \text{$i=t=\eta$, $j=s=\eta+1$},
\end{cases} \label{Engco}
\end{equation}
with the following restrictions on the $f_{ij}$s to ensure that all the
conditions of the twisting system are satisfied,
\begin{align}
f_{\eta\eta}&=f_{\eta+1,\eta+1}, \\
f_{\eta,\eta+1}&=q^{-1}p_{\eta,\eta+1}f_{\eta\eta}, \\
f_{\eta+1,\eta}&=q^{-1}p_{\eta+1,\eta}f_{\eta\eta}, \\
f_{i,\eta+1}&=p_{i,\eta+1}p_{\eta,i}f_{i\eta},\;\;\;\; i\neq\eta,\eta+1, \\
f_{\eta+1,i}&=p_{\eta+1,i}p_{i,\eta}f_{\eta i},\;\;\;\; i\neq\eta,\eta+1. \label{Engcon}
\end{align}
The new $R$-form is then given by the $R$-matrix
\begin{equation}
(R_{EK})_{ij}^{st}=\begin{cases}
q& \text{$i=j=s=t$},\\
\tilde{p}_{ij}& \text{$i=s\neq j=t$},\\
q-q^{-1}& \text{$i=t<j=s$}, \\
-(q-q^{-1})& \text{$i=t=\eta$, $j=s=\eta+1$}, \\
q-q^{-1}& \text{$i=t=\eta+1$, $j=s=\eta$},
\end{cases}
\end{equation}
where $\tilde{p}_{ij}=p_{ij}f_{ji}f_{ij}^{-1}$. There is no change in the number
of independent parameters. The twist quoted here actually corresponds to an 
embedding of $U_{q}(\mathfrak{gl}_{2}(\CC))$ in $U_{q}(\mathfrak{gl}_{n}(\CC))$.
There are many others, and we refer the reader to~\cite{Eng} for details.
\eex

\section{A new class of twisting $2$-cocycles}

Our major results all appear as particular examples of a new twisting system,
quite distinct from that of Reshetikhin, described in the following theorem.
\bth
Suppose $A(R)$ is any FRT bialgebra, defined in terms of an $n\times n$
$R$-matrix $R$. To any $n\times n$ matrix $F$ which satisfies the following
conditions,
\begin{align}
F_{12}F_{23}&=F_{23}F_{12}, \label{ncc1}\\
R_{12}F_{23}F_{13}&=F_{13}F_{23}R_{12}, \label{ncc2}\\
R_{23}F_{12}F_{13}&=F_{13}F_{12}R_{23}, \label{ncc3}
\end{align}
there corresponds a counital $2$-cocycle $\chi$ defined on $A(R)$. It is 
given on the generators of $A(R)$ by 
\begin{gather}
\chi(1,T)=\chi(T,1)=\ep(T), \label{ncc4}\\
\chi(T_1,T_2)=F_{12}, \label{ncc5}
\end{gather}
and extended to the whole algebra as
\begin{align}
\chi(T_{1}T_{2},T_{3})& =\chi(T_{2},T_{3})\chi(T_{1},T_{3})=F_{23}F_{13},  \label{ncc6} \\
\chi(T_{1},T_{2}T_{3})& =\chi(T_{1},T_{2})\chi(T_{1},T_{3})=F_{12}F_{13}. \label{ncc7}
\end{align}
\ethe
\begin{proof}
The fact that $\chi$ is consistent with the underlying algebraic structure of $A(R)$
follows from~(\ref{ncc2}) and~(\ref{ncc3}), while the defining condition~(\ref{coc2})
follows easily on using~(\ref{ncc1})
together with~(\ref{ncc6}),~(\ref{ncc7}) and the fact that 
$F_{ij}F_{kl}=F_{kl}F_{ij}$ whenever $i$, $j$, $k$ and $l$ are mutually distinct.
\end{proof}

\bre
There is of course a dual result to this, which applies
to any quasitriangular Hopf algebra $(H,\mathcal{R})$: Given an 
invertible element $\mathcal{F}\in H\ot H$ satisfying
\begin{equation}
\mathcal{F}_{12}\mathcal{F}_{23}=\mathcal{F}_{23}\mathcal{F}_{12},
\end{equation}
together with
\begin{align}
(\d\ot\id)(\mathcal{F})&=\mathcal{F}_{23}\mathcal{F}_{13},\\
\intertext{and}
(\id\ot\d)(\mathcal{F})&=\mathcal{F}_{12}\mathcal{F}_{13},
\end{align}
then $(H,\mathcal{R}_{\mathcal{F}})$ is a new quasitriangular Hopf algebra,
with the coproduct, universal $R$-matrix
and antipode twisted as in~(\ref{tc}),~(\ref{tr}) and~(\ref{ta})
respectively.
\ere

Some of the general features of twists coming from this construction will be
explicated in the following example.
\bex
Let us take as our initial object, the multiparameter standard quantum group 
$\CC_{q,p}[GL(3)]$, and
consider the possibility of defining on it a $2$-cocycle $\chi$ defined on the generators
as
\begin{equation}
\chi(T_{i}^{s},T_{j}^{t})=F_{ij}^{st}=\begin{cases}
f_{ij}& \text{$i=s$, $j=t$},\\
\mu& \text{$i=1$, $j=3$, $s=t=2$}.
\end{cases}
\end{equation}
For $F$ to satisfy~(\ref{ncc1}) we need $f_{i1}=f_{i2}$ and $f_{2i}=f_{3i}$ for all 
$i=1,\ldots,3$. For $\chi$ to be compatible with the algebra structure 
of $\CC_{q,p}[GL(3)]$ we need~(\ref{ncc2}) and~(\ref{ncc3}) to be satisfied,
which further requires $p_{i2}f_{i3}=p_{i1}f_{i2}$ and 
$p_{3i}f_{2i}=p_{2i}f_{1i}$ where $p_{ii}=q$, for $i=1,\ldots,3$. 
For generic $p_{ij}$ these equations have no solution. However,
giving up a degree of freedom from the parameter space of $\CC_{q,p}[GL(3)]$ by setting
$p_{13}=qp_{12}p_{23}$, they can be solved. As a matrix, $F$ is then given by
\begin{equation}
F=\begin{pmatrix}
q^{-1}p_{32}f & 0 & 0 & 0 & 0 & 0 & 0 & 0 & 0\\
0 & q^{-1}p_{32}f & 0 & 0 & 0 & 0 & 0 & 0 & 0\\
0 & 0 & p_{21}p_{32}f & 0 & \mu & 0 & 0 & 0 & 0\\
0 & 0 & 0 & f & 0 & 0 & 0 & 0 & 0\\
0 & 0 & 0 & 0 & f & 0 & 0 & 0 & 0\\
0 & 0 & 0 & 0 & 0 & q^{-1}p_{21}f & 0 & 0 & 0\\
0 & 0 & 0 & 0 & 0 & 0 & f & 0 & 0\\
0 & 0 & 0 & 0 & 0 & 0 & 0 & f & 0\\
0 & 0 & 0 & 0 & 0 & 0 & 0 & 0 & q^{-1}p_{21}f
\end{pmatrix},
\end{equation}
where $f=f_{22}$. The $R$-matrix of $\CC_{q,p}[SL(3)]$ with 
$p_{13}=qp_{12}p_{23}$ then twists to $R_{\chi}$, where
\begin{equation}
R_{\chi}=\begin{pmatrix}
q & 0 & 0 & 0 & 0 & 0 & 0 & 0 & 0\\
0 & qp & 0 & q-q^{-1} & 0 & 0 & 0 & 0 & 0 \\
0 & 0 & qp^{2} & 0 & -p^{2}qf^{-1}\mu & 0 & q-q^{-1} & 0 & 0 \\
0 & 0 & 0 & q^{-1}p^{-1} & 0 & 0 & 0 & 0 & 0 \\
0 & 0 & 0 & 0 & q & 0 & 0 & 0 & 0 \\
0 & 0 & 0 & 0 & 0 & qp & 0 & q-q^{-1} & 0 \\
0 & 0 & 0 & 0 & qf^{-1}\mu & 0 & q^{-1}p^{-2} & 0 & 0 \\
0 & 0 & 0 & 0 & 0 & 0 & 0 & q^{-1}p^{-1} & 0 \\
0 & 0 & 0 & 0 & 0 & 0 & 0 & 0 & q 
\end{pmatrix},
\end{equation}
with $p=p_{12}p_{23}$. It is clear that on choosing $f=-p\lambda^{-1}$ and
$\mu=q^{-1}(q-q^{-1})$, we have obtained precisely the $R$-matrix $R_{CG,p}$ 
for $n=3$.
\eex

The $2$-cocycle here is an example of a general class of
\emph{simple root $2$-cocycles} which are defined on $\CC_{q,p}[GL(n)]$ 
for any pair of integers $(k,l)$ such that $0<k<l<n$ by
\begin{equation}
\chi(T_{i}^{s},T_{j}^{t})=F_{ij}^{st}=\begin{cases}
f_{ij}& \text{$i=s$, $j=t$},\\
\mu& \text{$i=k$, $j=l+1$, $s=k+1$, $t=l$},
\end{cases}
\end{equation}
with the constraints
\begin{equation}
f_{i,k}=f_{i,k+1},\;\;\;\;\;\;f_{l,i}=f_{l+1,i},
\end{equation}
and
\begin{align}
p_{i,k}f_{i,l}&=p_{i,k+1}f_{i,l+1}, \\
p_{l,i}f_{k,i}&=p_{l+1,i}f_{k+1,i},
\end{align}
for all $i=1,\ldots,n$. The name comes from the fact that these twists add non-zero
elements to the $R$-matrix at points corresponding to the the non-zero elements
of the matrices $\Gamma(e_{\a_{k}})\otimes \Gamma(e_{-\a_{l}})$, 
and $\Gamma(e_{-\a_{l}})\otimes \Gamma(e_{\a_{k}})$, where 
the $e_{\a_{k}}$ are the basis elements corresponding to
the simple roots of the Lie algebra $\gln$, and $\Gamma$ is the first
fundamental representation. These simple root $2$-cocycles may be combined
in more general \emph{composite simple root $2$-cocycles} defined on 
$\CC_{q,p}[GL(n)]$ for each $0<k<n$, by
\begin{equation}
\chi(T_{i}^{s},T_{j}^{t})=F_{ij}^{st}=\begin{cases}
f_{ij}& \text{$i=s$, $j=t$},\\
\mu_{m}& \text{$i=k$, $j=m+1$, $s=k+1$, $t=m$}.
\end{cases}
\end{equation}
With the constraints as before, and $m$ now taking all possible values such that
$k<m<n$, this twist imparts a whole series of non-standard
off-diagonal elements to the $R$-matrix. Demonstration of the truth of these statements
involves a straightforward verification of the conditions~(\ref{ncc1}),~(\ref{ncc2})
and~(\ref{ncc3}).

While working on the universal
$T$-matrix, Fronsdal and Galindo~\cite{FG1,FG2} found  an interesting
non-standard deformation of $\CC[GL(2N-1)]$, which we shall
denote by $\CC_{FG,q}[GL(2N-1)]$, and
which they called `esoteric'. We have already introduced their $R$-matrix
in~(\ref{FGR}). For $N=2$ this is precisely the generalised Cremmer-Gervais
quantum group $\CC_{CG,q,p}[GL(3)]$ with $p=q^{-1}$ and
$\l=q^{2}(\kappa_{1}/(q-q^{-1}))$. However for $N>2$, their quantum groups 
do not coincide with those of the Cremmer-Gervais series (c.f.~added note 
in~\cite{FG2}). In fact, in a sense which
we will make more precise in the next section, the quantum groups
$\CC_{FG,q}[GL(2N-1)]$ are `not as non-standard' as those
of Cremmer and Gervais. As we shall discuss later, the general
Cremmer-Gervais quantum group does not seem to be a twisting of a
standard quantum group by a $2$-cocycle of the type we are considering
here. However the quantum groups of Fronsdal
and Galindo \emph{are} obtained from standard-type quantum groups through 
a $2$-cocycle which fits into our general scheme, and is presented in the
following proposition.
\bprop \label{fgcl}
On the quantum group $\CC_{q,p}[GL(2N-1)]$, with the parameters constrained
according to
\begin{equation}
p_{ji'}=qp_{jN}p_{Ni'},\;\;\;\;\;\frac{p_{ij}}{p_{iN}p_{Nj}}=
\frac{p_{i'j'}}{p_{i'N}p_{Nj'}}, \label{fgpc}
\end{equation}
for all $0<i,j<N$ and where $i'=2N-i$, there is a $2$-cocycle $\chi_{FG}$ defined as
\begin{equation}
\chi_{FG}(T_{i}^{s},T_{j}^{t})=F_{ij}^{st}=\begin{cases}
f_{ij}& \text{$i=s$, $j=t$},\\
\mu_{k}& \text{$i=k$, $j=k'$, $s=N$, $t=N$},\\
\lambda_{kl}& \text{$i=k$, $j=k'$, $s=l$, $t=l'$},
\end{cases} \label{fgc}
\end{equation}
where $0<k<l<N$.
All the $f_{ij}$ are given in terms of $f_{NN}$ according to
\begin{equation}
f_{ij}=\begin{cases}
q^{-1}p_{i'N}f_{NN}& \text{$0<i,j\leq N$},\\
p_{i'j}p_{jj'}f_{NN}& \text{$0<i\leq N<j<2N$},\\
f_{NN}& \text{$0<j\leq N<i<2N$},\\
q^{-1}p_{Nj'}f_{NN}& \text{$N<i,j<2N$}.\\
\end{cases} \label{fgfc}
\end{equation}
The $\lambda$s are determined in terms of the $\mu$s by
\begin{equation}
\lambda_{ij}=p_{j'j}f_{NN}(q-q^{-1})(\mu_{i}/\mu_{j}),
\end{equation}
for all $0<i<j<N$.
\eprop
\begin{proof}
This result is obtained by a applying conditions~(\ref{ncc1}),~(\ref{ncc2})
and~(\ref{ncc3}), with $R=R_{S,p}$, to~(\ref{fgc}).
\end{proof}
\bre
It is not difficult to establish that the number of parameters left
in $\CC_{q,p}[GL(2N-1)]$ after imposing the conditions~(\ref{fgpc}),
is $(1+\frac{1}{2}(N-1)(N+2))$. This is just the number of independent
$p_{ij}s$.
\ere

Using the $2$-cocycle~(\ref{fgc}) we in fact obtain an $R$-matrix more general
than that of Fronsdal and Galindo. Details are given in Appendix~B,
where we obtain this \emph{multiparameter generalised Fronsdal-Galindo $R$-matrix}
and demonstrate explicitly that their original
$R$-matrix~(\ref{FGR}) is a special case of the new generalised $R$-matrix.

\section{The Cremmer-Gervais problem for $GL(4)$ and beyond}
In the semiclassical theory of quasitriangular Lie bialgebras associated
with Lie algebras $\liea$, and 
their corresponding Poisson Lie groups (see for example the treatment in the 
book by Chari and Pressley~\cite{CP}), the fundamental role is played
by the classical $r$-matrix, $r\in\liea\ot\liea$, which completely specifies
the Lie bialgebra. In the case of complex, finite dimensional, simple
Lie algebras, there is a complete classification of all such $r$-matrices,
due to Belavin and Drinfeld~\cite{BD,CP}, in terms of `admissible' or 
`Belavin-Drinfeld' triples, $(\Pi_{1},\Pi_{0},\tau)$, where 
$\Pi$ is the set of simple roots of $\liea$, $\Pi_{1}$, $\Pi_{0}\subset\Pi$ and 
$\tau:\Pi_{1}\mapto\Pi_{0}$ is a bijection. Considering in particular the 
situation for $\liea=\sln$, we can distinguish three cases of
interest to us. In the standard, or Drinfeld-Jimbo case, the 
Belavin-Drinfeld triple has $\Pi_{1}$ and $\Pi_{0}$ both empty and the 
corresponding $r$-matrix, $r_{S}$, coincides with the semiclassical limit of the 
universal $R$-matrix, $\mathcal{R}_{S}$, of the familiar quasitriangular 
quantized universal enveloping algebra, $\qusln$. Another $r$-matrix, 
this time for $\slfg$, has~\cite{H4} 
$\Pi_{1}=\{\a_{1},\a_{2}, \ldots ,\a_{N-1}\}$ and 
$\Pi_{0}=\{\a_{N},\a_{N+1}, \ldots ,\a_{2(N-1)}\}$, where 
$\a_{1}, \ldots ,\a_{2(N-1)}$ are the simple roots of $\slfg$. When considered
in the first fundamental representation of $\slfg$, this can be seen to
correspond to the semiclassical limit of a Fronsdal-Galindo type 
$R$-matrix~(\ref{FGR}),
and so will be denoted $r_{FG}$. 
Finally, we have a $r$-matrix for $\sln$, in which $\Pi_{1}$ and
$\Pi_{0}$ are as full as possible~\cite{BDF}, with 
$\Pi_{1}=\{\a_{1},\a_{2}, \ldots ,\a_{n-2}\}$ and
$\Pi_{0}=\{\a_{2},\a_{n+1}, \ldots ,\a_{n-1}\}$, where 
$\a_{1}, \ldots ,\a_{n-1}$ are the simple roots of $\sln$.
This time, when viewed in the first fundamental representation of $\sln$,
we find a correspondence with a Cremmer-Gervais type $R$-matrix~(\ref{cg}), and
so we will write this $r$-matrix as $r_{CG}$. 
Let us note, that in each of these cases, the element $t$ defined as 
$t=r_{12}+r_{21}$ is identical, and is in fact the Casimir element of $\sln\ot\sln$.

In a series of fundamental works~\cite{D1,D3,D4}, Drinfeld proved 
that given any Lie algebra, $\liea$, together with a symmetric $\liea$-invariant
element, $t$, there exists a quantization of the universal enveloping algebra,
$U(\liea)$, as a \emph{quasitriangular quasi-Hopf quantized universal enveloping
algebra}, $(U(\liea)[[h]],\Phi,e^{ht/2})$, and that this quantization is 
\emph{unique up to twisting}. An immediate consequence of this result is that
the standard quantization, $(\qusln,\mathcal{R}_{S})$ of $U(\sln)$, is
twist equivalent, as a \emph{quasitriangular quasi-Hopf algebra}, to the 
`universal' quantization $(U(\sln)[[h]],\Phi,e^{ht/2})$. Subsequently~\cite{D2}, 
Drinfeld formulated a number of unsolved problems in quantum group theory. 
Among these, was the question of whether \emph{every} finite dimensional
Lie bialgebra admits a quantization as a quantized universal enveloping
algebra. This was recently answered, in the affirmative, 
by Etingof and Kazhdan~\cite{EK}.
Though their result did not provide an explicit construction, it does tell
us that in addition to the well known Drinfeld-Jimbo quasitriangular 
quantized universal enveloping algebra, $(\qusln,\mathcal{R}_{S})$, 
we must assume the \emph{existence} of
quasitriangular Fronsdal-Galindo and Cremmer-Gervais quantized universal
enveloping algebras, with corresponding universal $R$-matrices 
$\mathcal{R}_{FG}$ and $\mathcal{R}_{CG}$ respectively. Moreover, by Drinfeld's
earlier result, we know that these quantized universal
enveloping algebras must be twist equivalent as \emph{quasitriangular
Hopf algebras}, to $(\qusln,\mathcal{R}_{S})$. In particular, the universal 
$R$-matrices, $\mathcal{R}_{S}$, $\mathcal{R}_{FG}$, and $\mathcal{R}_{CG}$,
must each be related to each other by twisting in the style of 
equation~(\ref{tr}). 

It is reasonable, we believe, to work under the motivating assumption that the 
matrices $R_{FG}$ and $R_{CG}$, which have been considered in this paper,
correspond to the, as yet unknown, universal $R$-matrices, 
$\mathcal{R}_{FG}$ and $\mathcal{R}_{CG}$, in the first fundamental 
representation. In this case, the theory we have just outlined implies that 
in the dual world of co-quasitriangular Hopf algebras, there should exist 
$2$-cocycles for the construction of the Fronsdal-Galindo quantum groups 
and the Cremmer-Gervais quantum
groups from the standard quantum groups. Some support for the assumption has 
been provided in this paper, with the explicit construction of a twisting 
$2$-cocycle for the construction of the Fronsdal-Galindo quantum groups,
and the Cremmer-Gervais deformation of $GL(3)$. However the problem
for the Cremmer-Gervais deformations of $GL(n)$ for $n>3$ remains open.

The pair of non-standard off-diagonal elements which appear in the Cremmer-Gervais
$R$-matrix for $GL(3)$ `correspond' in the sense described above, to the
element $e_{\a_{1}}\wedge e_{-\a_{2}}$  of the corresponding classical 
$r$-matrix, $r_{CG}$. As we have seen, our twisting
construction has no problem dealing with this case. However, for $GL(4)$ and
beyond, the Cremmer-Gervais $R$-matrix involves an increasing number
of non-simple root combinations --- more than appear in the Fronsdal-Galindo
$R$-matrix, and our construction does not appear
to be able to deal with this circumstance. In particular, in the 
Cremmer-Gervais $R$-matrix for $GL(4)$, there are pairs of non-standard 
matrix elements corresponding to the $r$-matrix elements
$e_{\a_{1}}\wedge e_{-\a_{2}}$, $e_{\a_{1}}\wedge e_{-\a_{3}}$,
$e_{\a_{2}}\wedge e_{-\a_{3}}$ and $e_{\a_{1}+\a_{2}}\wedge e_{-(\a_{2}+\a_{3})}$.
The last term, in particular, causes problems. We finish by explaining how a \emph{new}
non-standard $GL(4)$ $R$-matrix may be obtained, which contains a pair of matrix
elements corresponding to this term.

Starting from the standard quantum group $\CC_{q,p}[GL(4)]$, 
and twisting first using a $2$-cocycle of the kind in~(\ref{Engco}), we
obtain a new quantum group given in terms of the $R$-matrix,
\begin{equation}
(R_{EK})_{ij}^{st}=\begin{cases}
q& \text{$i=j=s=t$},\\
\tilde{p}_{ij}& \text{$i=s\neq j=t$},\\
q-q^{-1}& \text{$i=t<j=s$}, \\
-(q-q^{-1})& \text{$i=t=2$, $j=s=3$}, \\
q-q^{-1}& \text{$i=t=3$, $j=s=2$}.
\end{cases}
\end{equation}
This quantum group is now amenable to a twist by one of our new $2$-cocycles,
given by
\begin{equation}
\chi(T_{i}^{s},T_{j}^{t})=F_{ij}^{st}=\begin{cases}
f_{ij}& \text{$i=s$, $j=t$},\\
\lambda& \text{$i=1$, $j=4$, $s=3$, $t=2$},
\end{cases}
\end{equation}
with the constraints,
\begin{gather}
f_{i1}=f_{i3},\;\;\;\;\;f_{2i}=f_{4i},\\
\tilde{p}_{i1}f_{i2}=\tilde{p}_{i3}f_{i4},\;\;\;\;\;
\tilde{p}_{4i}f_{3i}=\tilde{p}_{2i}f_{1i},
\end{gather}
for $i=1,\ldots,4$. Note that the this $2$-cocycle could not have been
defined on the original standard $R$-matrix. The $R$-matrix for this \emph{new}
non-standard quantum group may now be written as
\begin{equation}
(R_{NS})_{ij}^{st}=\begin{cases}
q& \text{$i=j=s=t$},\\
\gamma_{ij}& \text{$i=s\neq j=t$},\\
q-q^{-1}& \text{$i=t<j=s$}, \\
-(q-q^{-1})& \text{$i=t=2$, $j=s=3$}, \\
q-q^{-1}& \text{$i=t=3$, $j=s=2$},\\
\gamma_{14}\varrho& \text{$i=1$, $j=4$, $s=3$, $t=2$},\\
-q\gamma_{23}\varrho& \text{$i=4$, $j=1$, $s=2$, $t=3$},
\end{cases}
\end{equation}
where 
\begin{equation}
\gamma_{ij}=\tilde{p}_{ij}f_{ji}{f_{ij}}^{-1},\;\;\;\;\;
\varrho=-\lambda{f_{14}}^{-1}{f_{32}}^{-1},
\end{equation}
and 
\begin{equation}
\gamma_{12}\gamma_{23}=q\gamma_{24},\;\;\;\;\;\gamma_{24}\gamma_{34}=q\gamma_{14}.
\end{equation}
This new $R$-matrix depends on 6 parameters.
It might be interesting to investigate such double twists further.

\appendix
\section{The $3$-parameter generalised Cremmer-Gervais $R$-matrix}
\setcounter{equation}{0}
\renewcommand{\theequation}{\thesection\arabic{equation}}
We give here the derivation of the result quoted in~Example~\ref{mpcgex}.
We proceed in two stages, obtaining the required result by demonstrating that the
$R$-matrix~(\ref{mpcg}) is a twist of the $R$-matrix~(\ref{cg}).  

\noindent
{\bf{1. }} As explained in Section~2, as any diagonal matrix
$F_{ij}^{kl}=f_{ij}\delta_{i}^{k}\delta_{j}^{l}$ is a solution of the QYBE, we 
can define a twisting $2$-cocycle $\chi$ in terms of it as $\chi(T_{1},T_{2})=F_{12}$
as long as the compatibility conditions~(\ref{compat1}) and~(\ref{compat2}) are
satisfied. In terms of matrix components, these conditions become 
\begin{align}
R_{ij}^{st}f_{s\a}f_{t\a}&=f_{j\a}f_{i\a}R_{ij}^{st}, \\
\intertext{and}
R_{ij}^{st}f_{\a t}f_{\a s}&=f_{\a i}f_{\a j}R_{ij}^{st},
\end{align}
$i,j,\a,s,t=1,\ldots,n$, $n\geq 3$. Thus the non-zero elements of the $R$-matrix $R_{CG}$ 
determine the constraints on the elements of $F$. It is not difficult to see
that the only non-trivial relations which we get are
\begin{align}
f_{i\a}f_{j\a}&=f_{s\a}f_{t\a}\;\;\;\; i<s<j,\;\;t=i+j-s, \label{consist1}\\
\intertext{and}
f_{\a i}f_{\a j}&=f_{\a s}f_{\a t}\;\;\;\; i<s<j,\;\; t=i+j-s, \label{consist2}
\end{align}
$i,j,\a,s,t=1,\ldots,n$, $n\geq 3$. We will now prove the following lemma.
\ble
The system of equations~(\ref{consist1}) and~(\ref{consist2}), in $n^2$ unknowns,
has a solution space completely described in terms of four unknowns $x$, $y$, $z$, 
$w$ say, as
\begin{equation}
f_{ij}=x^{(i-2)(j-2)}y^{-(i-2)(j-1)}z^{-(i-1)(j-2)}w^{(i-1)(j-1)}, \label{tbp}
\end{equation}
$i,j=1,\ldots,n$.
\ele
\begin{proof}
We use induction. Consider the simplest case, $n=3$, so that we have $i=1$, 
$s=t=2$, $j=3$ and there are six equations
\begin{align}
f_{11}f_{31}&={f_{21}}^2, \label{sp1}\\
f_{12}f_{32}&={f_{22}}^2, \label{sp2}\\
f_{13}f_{33}&={f_{23}}^2, \label{sp3}\\
f_{11}f_{13}&={f_{12}}^2, \label{sp4}\\
f_{21}f_{23}&={f_{22}}^2, \label{sp5}\\
f_{31}f_{33}&={f_{32}}^2. \label{sp6}
\end{align}
Only five of these are independent, e.g.\ combining (\ref{sp1}), (\ref{sp2}),
(\ref{sp3}), (\ref{sp4}) and (\ref{sp5}) yields (\ref{sp6}), so the solution space
will be in terms of four unknowns. Choosing these to be $f_{11}=x$, $f_{12}=y$,
$f_{21}=z$ and $f_{22}=w$, we find
\begin{equation}
\|f_{ij}\|=
\begin{pmatrix} x&y&x^{-1}y^{2}\\z&w&z^{-1}w^{2}\\x^{-1}z^{2}&y^{-1}w^{2}&
xy^{-2}z^{-2}w^{4}\end{pmatrix},
\end{equation}
which verifies~(\ref{tbp}) for $n=3$. Now suppose that the solution space of the 
system of equations~(\ref{consist1}) and~(\ref{consist2}) for $n=k$, $k\geq 3$
is completely specified by~(\ref{tbp}), and consider $n=k+1$. Notice that we
still have all the equations we had for $n=k$ so~(\ref{tbp}) holds for
$i,j=1,\ldots,k$. We need to check that the new equations appearing 
consistently specify $f_{\a,k+1}$ and $f_{k+1,\a}$ according to~(\ref{tbp})
for $\a=1,\ldots,k+1$.

\noindent
From~(\ref{consist2}), for $\a=1,\ldots,k$,
\begin{align*}
f_{\a,k+1}&={f_{\a i}}^{-1}f_{\a s}f_{\a t},  \\ 
&=x^{-(\a-2)(i-2)+(\a-2)(s-2)+(\a-2)(t-2)} \\
&\times y^{(\a-2)(i-1)-(\a-2)(s-1)-(\a-2)(t-1)} \\
&\times z^{(\a-1)(i-2)-(\a-1)(s-2)-(\a-1)(t-2)} \\
&\times w^{-(\a-1)(i-1)+(\a-1)(s-1)+(\a-1)(t-1)}, \\ 
&=x^{(\a-2)(-i+s+t-2)}y^{(\a-2)(i-s-t+1)}z^{(\a-1)(i-s-t+2)}w^{(\a-1)(-i+s+t-1)}, \\ 
&=x^{(\a-2)(k+1-2)}y^{-(\a-2)(k+1-1)}z^{-(\a-1)(k+1-2)}w^{(\a-1)(k+1-1)},
\end{align*}
as required. Invoking the obvious symmetry between~(\ref{consist1}) 
and~(\ref{consist2}) we deduce the equivalent result from~(\ref{consist1})
for $f_{k+1,\a}$, $\a=1,\ldots,k$. Replacing $\a$ by $k+1$ in the above 
computation yields the correct result for $f_{k+1,k+1}$. The consistency of these
solutions still needs to be checked, but follows from the following. Take
$f_{\a,k+1}=f_{\a s}f_{\a t}{f_{\a i}}^{-1}$ from~(\ref{consist2}), and 
consider~(\ref{consist1}) with $\a=k+1$, i.e.\ 
\begin{equation*}
f_{i',k+1}f_{j',k+1}=f_{s',k+1}f_{t',k+1},
\end{equation*}
where $i'<s'<j'$, and $t'=i'+j'-s'$, $i',s',t',j'=1,\ldots,k+1$.
Then
\begin{align*}
\text{LHS }&=f_{i's}f_{i't}{f_{i'i}}^{-1}f_{j's}f_{j't}{f_{j'i}}^{-1} \\
&=f_{i's}f_{i't}f_{j'i}{f_{s'i}}^{-1}{f_{t'i}}^{-1}f_{j's}f_{j't}{f_{j'i}}^{-1} \\
&=f_{s's}f_{t's}f_{s't}f_{t't}{f_{s'i}}^{-1}{f_{t'i}}^{-1} \\
&=f_{s',k+1}f_{t',k+1} \\
&=\text{ RHS} 
\end{align*}
\end{proof}

\noindent
{\bf{2. }} We must now consider what combinations of the $f_{ij}$s actually
take part in the twisting. Thus, we consider the matrix 
$R_{CG,p}=F_{21}R_{CG}F^{-1}$, whose components are given by 
$(R_{CG,p})_{ij}^{st}=f_{ji}(R_{CG})_{ij}^{st}{f_{st}}^{-1}$. Explicitly,
\begin{equation}
(R_{CG,p})_{ij}^{st}=\begin{cases}
q& \text{$i=j=s=t$},\\
f_{ji}{f_{ij}}^{-1}qq^{-2(j-s)/n}& \text{$i=s<j=t$},\\
f_{ji}{f_{ij}}^{-1}q^{-1}q^{-2(j-s)/n}& \text{$i=s>j=t$}, \\
(q-q^{-1})& \text{$i=t<j=s$},\\
f_{ji}{f_{st}}^{-1}(q-q^{-1})q^{-2(j-s)/n}& \text{$i<s<j$, and $t=i+j-s$},\\
-f_{ji}{f_{st}}^{-1}(q-q^{-1})q^{-2(j-s)/n}& \text{$j<s<i$, and $t=i+j-s$}, 
\end{cases}
\end{equation}
and we are led to define
\begin{gather}
q_{ij}=f_{ij}{f_{ji}}^{-1}q^{-2(i-j)/n},\;\;\;\;\;i,j=1,\ldots,n, \\
\l_{ijst}=f_{ij}{f_{st}}^{-1}q^{-2(i-s)/n},\;\;\;\;\;i<s<j \text{ or } j<s<i \text{ and } t=i+j-s.
\end{gather}
The $q_{ij}$ and $\l_{ijst}$ satisfy the following obvious symmetries,
\begin{align}
q_{ji}&={q_{ij}}^{-1},\;\;\;\;\;i,j=1,\ldots,n, \\
\l_{jist}&=q_{ji}\l_{ijst},\;\;\;\;\;i<s<j \text{ and } t=i+j-s, \\
\l_{ijts}&=q_{st}\l_{ijst},\;\;\;\;\;i<s<j \text{ and } t=i+j-s.
\end{align}
Moreover, we see that 
modulo these symmetries every $\l_{ijst}$ must either be of the form
$\l_{ij\a \a}$, when $i+j$ is even, or $\l_{ij\a,\a+1}$ when
$i+j$ is odd \emph{or} be expressible in terms of these as
\begin{equation}
\l_{ijst}=\begin{cases}
\l_{ij\a\a}/\l_{st\a\a}& \text{$i+j$ even}, \\
\l_{ij\a,\a+1}/\l_{st\a,\a+1}& \text{$i+j$ odd}. \label{deco}
\end{cases}
\end{equation}
Now, recalling the solution~(\ref{tbp}), we find that that 
$q_{ij}=y^{-i+j}z^{i-j}q^{-2(i-j)/n}$,
so that on defining $p=y^{-1}zq^{-2/n}$, we have that $q_{ij}=p^{i-j}$. 
Now consider the
$\l_{ijst}$s. From~(\ref{tbp}) we get
\begin{align}
\l_{ij\a\a}&=
x^{-(\a-i)^{2}}y^{(\a-i)(\a-i+1)}z^{(\a-i)(\a-i-1)}w^{-(\a-i)^{2}}q^{-2(i-\a)/n} \nonumber\\
&=(y^{-1}zq^{-2/n})^{(i-\a)}(x^{-1}yzw^{-1})^{(\a-i)^{2}} \nonumber\\
&=p^{(i-\a)}(x^{-1}yzw^{-1})^{(\a-i)^{2}}, \\
\l_{ij\a,\a+1}&=
x^{-(\a-i)(\a-i+1)}y^{(\a-i)(\a-i+2)}z^{(\a-i)^{2}}w^{-(\a-i)(\a-i+1)}q^{-2(i-\a)/n} \nonumber\\
&=(y^{-1}zq^{-2/n})^{(i-\a)}(x^{-1}yzw^{-1})^{(\a-i)(\a-i+1)} \nonumber\\
&=p^{(i-\a)}(x^{-1}yzw^{-1})^{(\a-i)(\a-i+1)},
\end{align}
so that on defining $\l=x^{-1}yzw^{-1}$, and recalling~(\ref{deco}), we
find that
\begin{equation}
\l_{ijst}=p^{i-s}\l^{st-ij},\;\;\;\;\;i<s<j \text{ and } t=i+j-s.
\end{equation}
This completes the derivation.

\setcounter{equation}{0}
\renewcommand{\theequation}{\thesection\arabic{equation}}
\section{The Multiparameter Generalised Fronsdal-Galindo $R$-matrix}
In Proposition~\ref{fgcl} we presented the $2$-cocycle $\chi_{FG}$ defined on a 
certain standard quantum group in terms of a matrix $F$. To obtain the 
twisted quantum group we also need ${\chi_{FG}}^{-1}$, which is defined in terms
of $F^{-1}$,
\begin{equation}
(F^{-1})_{ij}^{st}=\begin{cases}
{f_{ij}}^{-1}& \text{$i=s$, $j=t$},\\
\bar{\mu}_{k}& \text{$i=k$, $j=k'$, $s=N$, $t=N$},\\
\bar{\lambda}_{kl}& \text{$i=k$, $j=k'$, $s=l$, $t=l'$},
\end{cases} 
\end{equation}
where $0<k<N$, $0<k<l<N$, and, 
\begin{align}
\bar{\mu}_{i}&=-qq^{i-i'}p_{ii'}{f_{NN}}^{-2}\mu_{i},\;\;\;\;\;0<i<N,\\
\bar{\lambda}_{ij}&=-q^{2(i-j)}p_{ii'}p_{jj'}{f_{NN}}^{-2}\lambda_{ij},
\;\;\;\;\;0<i<j<N.
\end{align}
Now we determine the twisted $R$-matrix, $R_{FG,p}$,
from $R_{FG,p}=F_{21}R_{S,p}F^{-1}$ as
\begin{equation}
(R_{FG,p})_{ij}^{st}=\begin{cases}
p_{ij}f_{ji}{f_{ij}}^{-1}& \text{$i=s$, $j=t$},\\
\bar{\mu}_{k}f_{k'k}p_{kk'}& \text{$i=k$, $j=k'$, $s=t=N$},\\
\bar{\lambda}_{kl}f_{k'k}p_{kk'}& \text{$i=k$, $j=k'$, $s=l$, $t=l'$},\\
\mu_{k}{f_{NN}}^{-1}p_{NN}& \text{$i=k'$, $j=k$, $s=t=N$},\\
\lambda_{kl}{f_{l'l}}^{-1}p_{l'l}& \text{$i=k'$, $j=k$, $s=l'$, $t=l$},\\
(q-q^{-1})& \text{$i=t<j=s$},
\end{cases} \label{mpfgr}
\end{equation}
where $0<k<N$ and $0<k<l<N$, the $p_{ij}$s and $f_{ij}$s are
constrained according to~(\ref{fgpc}) and~(\ref{fgfc}) respectively,
and all other parameters are given in terms of the $\mu$s.
We can refine the presentation of this $R$-matrix, setting 
\begin{align}
\k_{k}&=q^{-1}f_{NN}\bar{\mu}_{k},\\
\tilde{\k}_{k}&=-q^{2(N-k)}\k_{k},\\
\xi_{kl}&=(1-q^{2})(\k_{k}/\k_{l}),\\
\tilde{\xi}_{kl}&=(1-q^{-2})q^{2(l-k)}(\k_{k}/\k_{l}).
\end{align}
Then the $R$-matrix becomes the \emph{multiparameter
generalised Fronsdal-Galindo $R$-matrix},
\begin{equation}
(R_{FG,p})_{ij}^{st}=\begin{cases}
q& \text{$i=j=s=t$},\\
q{p_{ii'}}^{2}& \text{$i=s=j'$, $j=t$, $0<j<N$},\\
q^{-1}{p_{ii'}}^{2}& \text{$i=s$, $j=t=i'$, $0<i<N$},\\
p_{j'j}& \text{$i=s=N$, $j=t\neq N$},\\
p_{ii'}& \text{$i=s\neq N$, $j=t=N$},\\
p_{ij}p_{ii'}p_{j'i}& \text{$i=s\neq N$, $j=t\neq N$, $i\neq j$, $i+j\neq 2N$},\\
q-q^{-1}& \text{$i=t<j=s$},\\
qp_{ii'}\kappa_{i}& \text{$0<i<N$, $j=2N-i$, $s=t=N$},\\
qp_{j'j}\tilde{\kappa}_{j}& \text{$0<j<N$, $i=2N-j$, $s=t=N$},\\
q^{-1}p_{ii'}p_{ss'}\xi_{is}& \text{$0<i<s<N$, $j=2N-i$, $t=2N-s$},\\
qp_{j'j}p_{t't}\tilde{\xi}_{jt}& \text{$0<j<t<N$, $i=2N-j$, $s=2N-t$}.\\
\end{cases} 
\end{equation}
It is not difficult to see that, in general, this has $(1+\frac{1}{2}(N-1)(N+4))$
parameters. However, owing to the the particular $p_{ij}$s which 
appear in the $R$-matrix, in the cases of $N=2$ and $N=3$, the 
number of parameters is reduced to $3$ and $7$ respectively.

To identify the $R$-matrix originally discussed by Fronsdal
and Galindo~(\ref{FGR}), as a special case of this $R$-matrix, 
consider the particular solution 
of~(\ref{fgpc}) given by setting $p_{ij}=1$ for
$0<i\neq j<N$, $0<i<N<j<2N$ and $N<i\neq j<2N$.

\end{document}